\documentclass[reqno,11pt]{amsart}

\newtheorem{MainTheorem}{Theorem}
\newtheorem{Proposition}{Proposition}[section]
\newtheorem{Definition}[Proposition]{Definition}
\newtheorem{Lemma}[Proposition]{Lemma}
\newtheorem{Theorem}[Proposition]{Theorem}
\newtheorem{Corollary}[Proposition]{Corollary}
\newtheorem{Conjecture}[Proposition]{Conjecture}
\DeclareMathOperator{\vol}{vol}
\DeclareMathOperator{\svol}{svol}
\DeclareMathOperator{\Prob}{Prob}
\DeclareMathOperator{\Ext}{Ext}
\DeclareMathOperator{\conv}{conv}
\newcommand{\R}{\mathbb{R}}

\title{Centroid bodies and the convexity of area functionals}

\author{Andreas Bernig}
\email{bernig@math.uni-frankfurt.de}
\address{Institut f\"ur Mathematik, Goethe-Universit\"at Frankfurt, Robert-Mayer-Str. 10, 60054
Frankfurt, Germany}

\thanks{Supported by DFG grants BE 2484/3-1 and BE 2484/5-1.\\ AMS 2010 {\it Mathematics subject classification}:
52A38, 
52A40. 
}
\begin{document}

\begin{abstract}
We introduce a new volume definition on normed vector spaces. We show
that the induced $k$-area functionals are convex for all $k$. In the particular case $k=2$, our theorem implies that
Busemann's $2$-volume density is convex, which was recently shown by Burago-Ivanov. We also show how the new volume
definition
is related to the centroid body and prove some affine isoperimetric inequalities.  
\end{abstract}

\maketitle
\section{Introduction and statement of main results}

In a finite-dimensional Euclidean space, there is only one natural way to measure volumes of $k$-dimensional manifolds.
Similarly, there is basically only one natural volume definition on a Riemannian manifold. In contrast to this, measuring
volumes of submanifolds in a finite-dimensional normed space (or more generally on Finsler manifolds) is a more subtle
subject. Different aspects of the Euclidean volume give rise to different volume measurements. 

One natural way of defining volumes is to consider submanifolds in a normed space as metric spaces and to take the
corresponding Hausdorff measure. This gives rise to Busemann's definition of volume. Many basic questions, like
minimality of flat submanifolds, are still open. Recently, some progress was made by Burago and Ivanov
\cite{burago_ivanov12} who have shown that flat $2$-dimensional regions are minimal (see Corollary
\ref{cor_burago_ivanov} for the precise statement).  

A second well-known volume measurement is Holmes-Thompson volume, which equals the symplectic volume
of
the disc bundle. The use of symplectic geometry gives rise to a number of interesting results. It was shown recently by Ludwig \cite{ludwig10} that the
Holmes-Thompson surface area can be uniquely characterized by a valuation property. However, Holmes-Thompson
volume lacks some basic convexity properties. 

In geometric measure theory it is common to use Gromov's mass*, which has very strong convexity properties, but seems
less natural from the point of view of convex geometry. 

In this paper, we propose a new natural definition of volume which is based on a version of the well-known
centroid body and was inspired by a recent result of Burago-Ivanov \cite{burago_ivanov12}. We show that our definition
of
volume has strong convexity properties. More precisely, it induces convex $k$-densities for all $k$. Since the
$2$-volume density induced by our definition of volume equals the Busemann $2$-volume density, we obtain as a corollary 
Burago-Ivanov's theorem that $2$-planes are minimal with respect to Busemann volume. 

Let us describe our results in more detail. References for this section are \cite{alvarez_thompson} and
\cite{thompson_book96}. We let $\Lambda^k V$ denote the $k$-th exterior power of $V$ and $\Lambda^k_s V$ the cone of
simple $k$-vectors. 

\begin{Definition}[Definition of volume]
A definition of volume $\mu$ assigns to each normed vector space $(V,\|\cdot\|)$ a norm $\mu_V$ on
$\Lambda^nV$ (where $n=\dim V$) such that the following two conditions are satisfied:
\begin{enumerate}
 \item If $V$ is Euclidean, then $\mu_V$ is induced by the usual Lebesgue measure.
\item If $f:(V,\|\cdot\|) \to (W,\|\cdot\|)$ is a linear map that does not increase distances, then the induced map
$\Lambda^nf:(\Lambda^nV,\mu_V) \to (\Lambda^nW,\mu_W)$ does not increase distances. 
\end{enumerate}
\end{Definition}

If we want to stress the dependence on the norm, we will write $\mu_B$ instead of $\mu_V$, where $B$ is the unit ball in $(V,\|\cdot\|)$. 

From i) and ii) it follows that the map $(V,\|\cdot\|) \mapsto (\Lambda^n V, \mu_V)$ is continuous with respect to
the Banach-Mazur distance. 

An equivalent definition is as follows. Let $\mathcal{K}_0^s$ be the space of centrally symmetric compact convex bodies with non-empty interior. 

Given a definition of volume definition $\mu$ on $V$, define the functional $\mathcal{V}:\mathcal{K}_0^s \to \R_+$ by  
\begin{displaymath}
 \mathcal{V}(B):=\mu_B(B).
\end{displaymath}
The functional $\mathcal{V}$ satisfies the following properties:
\begin{enumerate}
 \item $\mathcal{V}$ is invariant under linear maps, i.e. $ \mathcal{V}(gB)= \mathcal{V}(B)$ for all $g \in
\mathrm{GL}(V)$. 
\item $ \mathcal{V}(E)=\omega_n$ (the usual volume of the Euclidean unit ball) if $E$ is an ellipsoid. 
\item \label{item_defvol_mon} If $B \subset B'$, then 
\begin{displaymath}
 \frac{ \mathcal{V}(B)}{\vol B} \geq \frac{ \mathcal{V}(B')}{\vol B'}
\end{displaymath}
where $\vol$ denotes any Lebesgue measure on $V$.
\end{enumerate}
Conversely, any functional with these properties defines a definition of volume. 

We call $ \mathcal{V}$ the {\it associated affine invariant}. 

\begin{Definition}[Main examples of definitions of volume]
\begin{enumerate}
 \item The Busemann definition of volume \cite{busemann47} has the associated affine invariant 
\begin{displaymath}
  \mathcal{V}^b(B) \equiv \omega_n.
\end{displaymath}
It equals the Hausdorff measure of $B$ with respect to the metric induced by $B$.  
\item Holmes-Thompson definition of volume \cite{holmes_thompson79} has associated affine invariant
\begin{displaymath}
  \mathcal{V}^{ht}(B):=\frac{1}{\omega_n} \svol(B \times B^\circ).
\end{displaymath}
Here $\svol$ denotes the symplectic volume on $V \times V^*$ and $B^\circ \subset V^*$ is the polar body of $B$. 
\item The Benson definition of volume (also called Gromov's mass*, see \cite{benson62, gromov83}) has associated affine invariant 
\begin{displaymath}
 \mathcal{V}^{m*}(B):= 2^n \frac{\vol B}{\inf_{P \supset B} \vol P} 
\end{displaymath}
Here $P$ ranges over all parallelotopes circumscribed to $B$ and $\vol$ is any choice of
Lebesgue measure on
$V$.
\item Ivanov's definition of volume \cite{ivanov08} has associated affine invariant
\begin{displaymath}
\mathcal{V}^i(B):= \omega_n \frac{\vol B}{\vol E},
\end{displaymath}
where $E$ is the maximal volume ellipsoid inscribed in $B$ (i.e. the John ellipsoid). 
\end{enumerate}
\end{Definition}

Each definition of volume on $V$ induces {\bf $k$-volume densities}, i.e. $1$-homogeneous, continuous, positive
functions on the set
of simple $k$-vectors in $V$, where $0 \leq k \leq \dim V$. More precisely, given a simple $k$-vector $a$, we put
\begin{displaymath}
 \mu_k(a):=\mu_{\langle a \rangle}(a),
\end{displaymath}
where $\langle a \rangle$ is the $k$-dimensional space spanned by $a$ with the induced norm. 

\begin{Definition}
A $k$-volume density $\mu_k$ is called {\it
extendibly convex} if it is the restriction of a norm on $\Lambda^k V$. 
\end{Definition}

There are other notions of convexity for $k$-volume densities. The $k$-volume density $\mu_k$ is called {\it totally
convex} if for each $k$-subspace in $V$, there exists a $\mu_k$-decreasing linear projection
onto that subspace. It is called {\it semi-elliptic}, if a plane $k$-disc has minimal $\mu_k$-area among all
Lipschitz
chains with the same boundary. Semi-ellipticity depends in a subtle way on the choice of the coefficient
ring. Semi-ellipticity over $\R$ is equivalent to extendible convexity \cite{burago_ivanov04}. This notion is important
in geometric measure theory, in particular in the solution of the Plateau problem in normed or metric spaces
\cite{almgren68, ambrosio_kirchheim00, gromov83, wenger08}. 

In codimension $1$, these notions coincide. In general, total convexity implies extendible
convexity. We refer to \cite{alvarez_thompson} for more
details and other notions of convexity. 

The $k$-density induced by the Holmes-Thompson volume is extendibly convex for $k=\dim V-1$. Busemann-Ewald-Shephard \cite{busemann_ewald_shephard} and later
Burago and Ivanov \cite{burago_ivanov02} gave examples showing that for $1<k<n-1$, it is not necessarily extendibly convex.

Busemann's volume is also convex in codimension $1$, but this is more difficult to show (in fact this is equivalent to
Busemann's intersection theorem \cite{busemann49}). 
An open conjecture (which appears as problem number 10 on Busemann-Petty's list \cite{busemann_petty56} of
problems in convex geometry) states that Busemann's definition of volume induces totally convex $k$-volume densities
for all $k$. The case $k=2$ of this
conjecture
was recently confirmed by Burago and Ivanov \cite{burago_ivanov12}. 

Gromov's mass* and Ivanov's definition of volume have the best convexity properties, as the induced $k$-volume densities are
totally convex for all $k$
\cite{benson62, gromov83, ivanov08}. 

For closely related results on the minimality of totally geodesic submanifolds of Finsler manifolds we refer to
\cite{alvarez_berck, berck09b, ivanov01, schneider01}. 

The aim of the present paper is to introduce a new definition of volume which has strong convexity properties. For
simplicity, we will use a fixed Lebesgue measure $\vol$ on $V$. Let 
$V(K_1,\ldots,K_n)$ denote the associated mixed volume of the compact convex bodies $K_1,\ldots,K_n$. We will follow the
usual
notation and write
\begin{displaymath}
 V_i(K,L):=V(K[n-i],L[i])=V(\underbrace{K,\ldots,K}_{n-i},\underbrace{L,\ldots,L}_i).
\end{displaymath}
The projection body of a compact convex body $K \subset V$ will be denoted by $\Pi K \subset V^*$, see the next section
for the definition and some properties.

Our main theorem is the following.

\begin{MainTheorem} \label{mainthm_convexity}
Let 
\begin{equation} \label{eq_def_functional}
  \mathcal{V}(B):=\frac{\omega_n^{n-1}}{\omega_{n-1}^n} \sup_{L \in \mathcal{K}(V)} \left\{\svol(B \times \Pi L)\left|
V_1(L,B)=1\right.\right\}.
\end{equation}
Then $ \mathcal{V}$ is the associated affine invariant of a definition of volume. The induced $k$-volume densities are
extendibly convex for all $k$. Moreover, 
\begin{displaymath}
  \mathcal{V}(B) \geq  \mathcal{V}^{ht}(B)
\end{displaymath}
for all unit balls $B$, with equality precisely for ellipsoids. 
\end{MainTheorem}

To the best of our knowledge, only two other definitions of volume with extendibly convex densities were known
previously, namely Gromov's mass* and Ivanov's definition of volume.

Theorem \ref{mainthm_convexity} implies a recent result by Burago
and Ivanov. 

\begin{Corollary}[Burago-Ivanov, \cite{burago_ivanov12}] \label{cor_burago_ivanov}
 The $2$-volume density induced by Busemann's definition of volume is extendibly convex.
\end{Corollary}

We will also show (Proposition \ref{prop_equivalent_conj}) that the stronger inequality $\mathcal{V}(B) \geq
\mathcal{V}^b(B)$ is equivalent to Petty's conjectured projection inequality (Conjecture \ref{conj_petty}). 

Our second main theorem establishes a link between our new definition of volume, the centroid body and random simplices. 

Recall that the support function of a convex body $K \subset V$ is the function $h(K,\cdot): V^* \to \R, \xi \mapsto \sup_{x \in K} \langle \xi,x\rangle$. 

\begin{MainTheorem}[Alternative description of $ \mathcal{V}$] \label{mainthm_gamma_mu}
Let $B \subset V$ be the unit ball of some norm. Let $\nu$ be a probability measure on $B^\circ$. Define a convex body
$\Gamma_\nu B^\circ \subset V^*$ by 
\begin{displaymath}
 h(\Gamma_\nu B^\circ,u):=\int_{B^\circ} |\langle \xi,u\rangle|d\nu(\xi), \quad u \in V.
\end{displaymath} 
\begin{enumerate}
 \item \label{item_omega_and_zmu} We have 
\begin{displaymath}
  \mathcal{V}(B) = \frac{\omega_n^{n-1}}{\omega_{n-1}^n} \left(\frac{n}{2}\right)^n \sup_{\nu \in \Prob(B^\circ)}
\svol(B \times \Gamma_\nu B^\circ).
\end{displaymath}
\item Let $[0,\xi_1,\ldots,\xi_n]$ be the simplex spanned by
$\xi_1,\ldots,\xi_n \in V^*$. Then 
\begin{displaymath}
 \vol(\Gamma_\nu B^\circ)= 2^n \int_{B^\circ} \cdots \int_{B^\circ} \vol [0,\xi_1,\ldots,\xi_n] d\nu(\xi_1) \ldots
d\nu(\xi_n).
\end{displaymath}
\item There exists a unique even probability measure on $B^\circ$ which maximizes $\vol \Gamma_\nu B^\circ$. It is
supported
in the set $\Ext B^\circ$ of extremal points of $B^\circ$. 
\end{enumerate}
\end{MainTheorem}

Remark: if $\nu$ is the uniform measure on $B^\circ$, then $\Gamma_\nu B^\circ=\Gamma B^\circ$, the well-known centroid
body \cite{gardner_book06, lutwak86, lutwak90}. In general, we call $\Gamma_\nu B^\circ$ the {\it centroid
body with respect to $\nu$}. 

\subsubsection*{Acknowledgements}
Some parts of this paper were worked out during a stay at the Universit\'e de Fribourg and I thank Stefan Wenger for
very
fruitful discussions. I also would like to thank Monika Ludwig, Rolf Schneider, Franz Schuster, Deane Yang and the
anonymous referee for useful remarks. 

\section{Notations and background}

We refer to the books by Schneider
\cite{schneider_book93} and Gardner \cite{gardner_book06} for information on convexity and the Brunn-Minkowski theory.
Let us recall some notions and theorems which will be used later on. 
 
Let $V$ be a real vector space of dimension $n$. The space of compact convex bodies in
$V$ is denoted by $\mathcal{K}(V)$. The space of symmetric compact convex bodies with non-empty interior will be denoted
by $\mathcal{K}_0^s(V)$. The convex hull of a set $X \subset V$ will be denoted by $\conv X$. 

A set $B \subset
\mathcal{K}_0^s(V)$ is the unit ball of some norm on $V$ and vice versa. If
$B=E$ is an ellipsoid, then the corresponding norm is Euclidean. 

Each $K \in \mathcal{K}(V)$ may be described by its support function $h(K,\xi):=\sup_{x \in K} \xi(x), \xi \in V^*$.
For $B \in \mathcal{K}_0^s$, the radial function is defined by $\rho(B,v):=\sup\{ \lambda \geq 0, \lambda v \in B\}, v
\in V, v \neq 0$. Note that $h$ is $1$-homogeneous, while $\rho$ is $(-1)$-homogeneous. 

The {\it polar body} of $B \in \mathcal{K}_0^s(V)$ is defined by 
\begin{displaymath}
 B^\circ:=\left\{\xi \in V^*: \xi(x) \leq 1, \forall x \in B\right\} \subset V^*. 
\end{displaymath}
We have 
\begin{displaymath}
 h(B,\xi)=\frac{1}{\rho(B^\circ,\xi)}, \xi \in V^*, \xi \neq 0. 
\end{displaymath}

The {\it mixed volume} of compact convex bodies will be denoted by $V(K_1,\ldots,K_n)$ and we will abbreviate $V_i(K,L):=V(K[n-i],L[i])$. We will use the following
inequality of Minkowski, which is a special case of the Alexandrov-Fenchel inequality:
\begin{displaymath}
 V_1(K,L)^n \geq \vol(K)^{n-1}\vol(L) .
\end{displaymath}
If $K,L$ contain inner points, then equality holds if and only if $K$ and $L$ are homothetic. 

Let us also recall the Brunn-Minkowski inequality:
\begin{displaymath}
 \vol(\lambda K +(1-\lambda)L)^\frac{1}{n} \geq \lambda \vol(K)^\frac{1}{n}+(1-\lambda) \vol(L)^\frac{1}{n}, \quad 0
\leq \lambda \leq 1. 
\end{displaymath}
If $K$ and $L$ contain interior points and $0<\lambda<1$, then equality holds if and only if $K$ and $L$ are
homothetic. 

The space $V \times V^*$ admits a symplectic volume form $\svol$ \cite{cannas_da_silva_book}. If $\vol$ is any choice of
Lebesgue measure on $V$ and
$\vol^*$ the dual Lebesgue measure on $V^*$, then $\svol=\vol \times \vol^*$. 

Let $V$ be a vector space, $\Omega \in \Lambda^n V^*$ a volume form and $K \in \mathcal{K}(V)$. The
{\it projection
body} \cite{gardner_book06} $\Pi K \in \mathcal{K}(V^*)$ is defined as follows. Let $v \in V, v \neq 0$. Then
$i_v\Omega:=\Omega(v,\cdot)$ is
a volume form on
the quotient space $V/\R \cdot v$ and $h(\Pi K,v):=\vol(\pi_v K,i_v\Omega)$, where $\pi_v:V \to V/\R \cdot v$ is
the quotient map, defines the support function of $\Pi K$. We will write $\Pi^\circ K:=(\Pi K)^\circ$ for
the polar body of $\Pi K$.  

Let us recall a well-known geometric inequality related to the projection body. 

\begin{Theorem}[Petty's projection inequality, \cite{petty71}] \label{thm_ppi}
Let $K \subset V$ be a compact convex body and $E \subset V$ an ellipsoid. Then 
\begin{displaymath}
 \vol(K)^{n-1} \vol \Pi^\circ K \leq \vol(E)^{n-1} \vol \Pi^\circ E 
\end{displaymath}
with equality precisely for ellipsoids. 
\end{Theorem}

The following conjecture is a strengthening of Petty's projection inequality. We refer to \cite{campi_gronchi06,
lutwak90b, martini_mustafaev10, petty71} for more information and equivalent formulations. 

\begin{Conjecture}[Petty's conjectured projection inequality] \label{conj_petty}
Let $K$ be a compact convex body in $V$ and $E \subset V$ an ellipsoid. Then 
\begin{displaymath}
 \vol(\Pi K) \vol (K)^{1-n} \geq \vol(\Pi E) \vol (E)^{1-n}
\end{displaymath}
with equality precisely for ellipsoids.  
\end{Conjecture}

The {\it centroid body} $\Gamma K \in \mathcal{K}(V)$ of a compact convex $K$ with non-empty interior is defined by 
\begin{displaymath}
 h(\Gamma K,\xi)=\frac{1}{\vol K} \int_K |\langle \xi,u\rangle| du, \quad \xi \in V^*.
\end{displaymath}
It satisfies the Busemann-Petty centroid inequality \cite{gardner_book06, lutwak_handbook93, schuster07}
\begin{equation} \label{eq_busemann_petty_centroid}
 \vol(\Gamma K) \geq \left(\frac{2\omega_{n-1}}{(n+1)\omega_n}\right)^n \vol K,
\end{equation}
with equality precisely for ellipsoids centered at the origin. 
 
Let $V$ be a Euclidean vector space with unit sphere $S^{n-1}$. The {\it cosine transform} is defined by 
\begin{displaymath}
 Cf(v):=\int_{S^{n-1}} |\langle u,v\rangle| f(u) d\sigma(u), \quad f \in C(S^{n-1}).
\end{displaymath}
where $d\sigma$ denotes the spherical Lebesgue measure. On the space of even, smooth functions, the cosine
transform is a bijection. The
cosine transform extends to measures on the sphere by 
\begin{displaymath}
 C \nu(v):=\int_{S^{n-1}} |\langle u,v\rangle| d\nu(u). 
\end{displaymath}
The cosine transform is injective on the space of even measures on $S^{n-1}$ (\cite{gardner_book06}, Appendix C.2). 

The space of probability measures on a topological space $X$ will be denoted by $\Prob(X)$.

\section{Proof of Theorem \ref{mainthm_gamma_mu}}

Let $B \subset V$ be the unit ball of some norm. By (\cite{schneider_book93}, 5.3.38), 
\begin{equation} \label{eq_functional_as_volume}
\vol(\Gamma_\nu B^\circ)=2^n \int_{B^\circ} \cdots \int_{B^\circ} \vol[0,\xi_1,\ldots,\xi_n] d\nu(\xi_1) \ldots
d\nu(\xi_n).
\end{equation}

\begin{Proposition} \label{prop_uniqueness}
There exists a unique
even probability measure $\nu$ on $B^\circ$ which maximizes $\vol(\Gamma_\nu B^\circ)$.
\end{Proposition}

\proof
By Prokhorov's theorem (see e.g. \cite{billingsley_book99}, Thm. 5.1), the space of probability measures on $B^\circ$ is
sequentially compact with
respect to weak convergence. Since the functional $\nu \mapsto \vol(\Gamma_\nu B^\circ)$ is continuous with respect to
weak topology, it follows that the supremum is attained. 

If the measure of the interior of $B^\circ$ is positive, then radial projection from $B^\circ$ onto $\partial
B^\circ$ (with the origin mapped to an arbitrary boundary point) of $\nu$ will increase our functional, hence each
optimal
measure $\nu$ must be concentrated on the boundary. Moreover, replacing $\nu$ by
its even part $\nu^{ev}$ does not change $\Gamma_\nu$, hence we may assume that $\nu$ is even, i.e. invariant under
central symmetry.  

Let $\nu,\tau \in \Prob(\partial B^\circ)$ be even measures and $0 < \lambda < 1$. Then
$\Gamma_{(1-\lambda)\nu+\lambda
\tau}B^\circ=(1-\lambda)\Gamma_\nu B^\circ+\lambda
\Gamma_\tau B^\circ$. By the Brunn-Minkowski inequality, it follows that 
\begin{displaymath}
\left( \vol \Gamma_{(1-\lambda)\nu+\lambda \tau}B^\circ\right)^\frac{1}{n} \geq (1-\lambda) \vol(\Gamma_\nu
B^\circ)^\frac{1}{n}+\lambda
\vol(\Gamma_\tau B^\circ)^\frac{1}{n}.
\end{displaymath}

This shows that the function $\nu \mapsto \vol(\Gamma_\nu B^\circ)^\frac{1}{n}$ is concave on the space of even
measures on $\partial B^\circ$. If $\vol \Gamma_\nu
B^\circ>0$, then equality in the
above inequality holds if and only if $\Gamma_\nu B^\circ$ is homothetic to $\Gamma_\tau B^\circ$, i.e. $\Gamma_\nu
B^\circ=t \Gamma_\tau B^\circ+v$ for $t>0, v \in V$. 

We claim that this is possible only if $\nu=\tau$. Indeed, since $\Gamma_\nu B^\circ$ and $\Gamma_\tau B^\circ$ are
centrally symmetric, $v=0$. Choose a Euclidean scalar product
on $V^*$ with unit sphere $S^{n-1}$. Let $\tilde \nu, \tilde \tau$ be the push-forwards of $\nu$ and $\tau$ under the
radial projection $\partial B^\circ \to S^{n-1}$. From $\Gamma_\nu B^\circ=t \Gamma_\tau B^\circ$ and from the injectivity of the 
cosine transform on even measures we deduce that $d\tilde \nu=t d \tilde \tau$ and hence $\nu=t \tau$. Since
$\nu$ and $\tau$ are probability measures, $t=1$.  

From this the uniqueness of the maximum follows easily. 
\endproof

Recall that a point $x$ in a compact convex body $K$ is called an {\it extreme point} if it cannot be written as
$x=\frac{a+b}{2}$ with $a,b \in K, a \neq b$. The set of extremal points is denoted by $\Ext K$. We refer to
\cite{barvinok} for more information. 
 
\begin{Proposition} \label{prop_support_mu}
 Let $B$ be a unit ball. The even measure $\nu$ such that $\vol \Gamma_\nu B^\circ$ is maximal is concentrated in the
set $\Ext B^\circ$ of extremal points.
\end{Proposition}

\proof
Let
\begin{displaymath}
\Delta:=\left\{(\lambda_1,\ldots,\lambda_{n+1})\in \R^{n+1}| \lambda_i \geq 0, \sum_i \lambda_i=1\right\} 
\end{displaymath}
be the standard simplex. By Minkowski's theorem (\cite{barvinok}, Thm. II.3.3), $B^\circ=\mathrm{conv} \Ext B^\circ$. Carath\'eodory's theorem (\cite{barvinok}, Thm. I.2.3) implies
that the
continuous map
\begin{align*}
 \left(\Ext B^\circ\right)^{n+1} \times \Delta & \to B^\circ \\
(\eta_1,\ldots,\eta_{n+1},\lambda_1,\ldots,\lambda_{n+1}) & \mapsto \sum_i \lambda_i \eta_i
\end{align*}
is onto. 

We fix a measurable right inverse $\xi \mapsto 
(\eta_1(\xi),\ldots,\eta_{n+1}(\xi),\lambda_1(\xi),\ldots,\lambda_{n+1}(\xi))$ of this map. Then each $\lambda_i:B^\circ \to
\R$ is a non-negative measurable function and each $\eta_i:B^\circ \to \Ext B^\circ$ is a measurable map. 
 
Let $\nu_i:=\nu \llcorner \lambda_i$ be the measure on $B^\circ$ with density function $\lambda_i$ with respect to $\nu$. Define a probability measure $\tilde \nu$ on $\Ext B^\circ \subset B^\circ$ by 
\begin{displaymath}
\tilde \nu:= \sum_{i=1}^{n+1} (\eta_i)_* \nu_i,
\end{displaymath}
where $(\eta_i)_*$ denotes the push-forward.   

If $f:B^\circ \to \R$ is a convex function, then 
\begin{align*}
 \int_{B^\circ} f(\xi) d\nu(\xi) & = \int_{B^\circ} f\left(\sum_i \lambda_i(\xi) \eta_i(\xi)\right) d\nu(\xi)\\
& \leq \int_{B^\circ} \sum_i \lambda_i(\xi) f(\eta_i(\xi)) d\nu(\xi)\\
& = \int_{B^\circ} \sum_i f(\eta_i(\xi)) d\nu_i(\xi)\\
& = \int_{B^\circ} f(\xi) d\tilde \nu(\xi).
\end{align*}

Since the function $\vol [0,\xi_1,\ldots,\xi_n]$ is convex in each
variable $\xi_i$, it follows from \eqref{eq_functional_as_volume} that 
\begin{displaymath}
\vol \Gamma_\nu B^\circ \leq \vol \Gamma_{\tilde \nu} B^\circ.
\end{displaymath}
By the uniqueness of the optimal measure, $\nu$ equals the even part of $\tilde \nu$ and is therefore 
concentrated on $\Ext B^\circ$. 
\endproof

\begin{Proposition}
 We have 
\begin{displaymath}
 \sup_{\nu \in \Prob(B^\circ)} \vol \Gamma_\nu B^\circ=\left(\frac{2}{n}\right)^n \sup\left\{\vol \Pi L: L \in
\mathcal{K}(V), V_1(L,B)=1\right\}.
\end{displaymath}
\end{Proposition}

\proof
We will use an auxiliary scalar product on $V$, which allows us to identify $V$ and $V^*$. Let $\nu \in \Prob(B^\circ)$ maximize $\vol \Gamma_\nu B^\circ$. We may assume that $\nu$ is concentrated on the boundary of $B^\circ$ and even.  

Let $\tilde \nu$ be the
push-forward of $\nu$ under the radial projection $\partial K^\circ \to S^{n-1}$. 

By the solution of Minkowski's problem applied to $\tilde \nu$, $\Gamma_\nu B^\circ$ is a projection body, say $\Gamma_\nu B^\circ=\Pi \tilde L$ for
some centrally symmetric convex body $\tilde L \subset V$.  

Then for $u \in V$
\begin{align*}
h(\Gamma_\nu B^\circ,u) & =  \int_{\partial B^\circ} |\langle \xi,u\rangle|d \nu(\xi)\\
& = \int_{S^{n-1}} \rho(B^\circ,\xi) |\langle \xi,u\rangle| d\tilde \nu(\xi).
\end{align*}

On the other hand, 
\begin{displaymath}
 h(\Pi \tilde L,u) = \frac12 \int_{S^{n-1}} |\langle \xi,u\rangle| dS_{n-1}(\tilde L,\xi),
\end{displaymath}
where $S_{n-1}(\tilde L,\cdot)$ is the surface area measure of $\tilde L$.

Using the injectivity of the cosine transform on even measures on the sphere, we find 
\begin{displaymath}
 d\tilde \nu=\frac12 \rho(B^\circ,\cdot)^{-1} \cdot
dS_{n-1}(\tilde L,\cdot)=\frac12 h(B,\cdot) dS_{n-1}(\tilde L,\cdot).
\end{displaymath}
Since $\tilde \nu$ is a probability measure, we must have 
\begin{displaymath}
 1=\int_{S^{n-1}} d\tilde \nu(\xi) = \frac12 \int_{S^{n-1}} h(B,\xi) dS_{n-1}(\tilde L,\xi)=\frac{n}{2}
V_1(\tilde L,B).   
\end{displaymath}

Let $L:=\left(\frac{n}{2}\right)^\frac{1}{n-1} \tilde L$. Then $V_1(L,B)=1$ and $\vol \Gamma_\nu
B^\circ=\vol \Pi \tilde L=\left(\frac{2}{n}\right)^n \vol
\Pi L$. Thus we have the inequality
\begin{displaymath}
 \sup_{\nu \in \Prob(B^\circ)} \vol \Gamma_\nu B^\circ \leq \left(\frac{2}{n}\right)^n \sup\left\{\vol \Pi L: L \in
\mathcal{K}(V), V_1(L,B)=1\right\}.
\end{displaymath}

To prove the inverse inequality, take $L$ with $V_1(L,B)=1$ and set $\tilde
L:=\left(\frac{2}{n}\right)^\frac{1}{n-1} L$. We define $d \tilde \nu:=\frac12
h(B,\cdot)
dS(\tilde L,\cdot)$, which is a probability measure on $S^{n-1}$. If $\nu$ is the push-forward of $\tilde \nu$ under the
radial
projection
$S^{n-1} \to \partial B^\circ$, then $\Gamma_\nu B^\circ=\Pi \tilde L$ and $\vol \Gamma_\nu
B^\circ=\left(\frac{2}{n}\right)^n \vol \Pi L$. 
\endproof


\section{Proof of Theorem \ref{mainthm_convexity}}

\begin{Lemma}
 The functional $ \mathcal{V}$ defined by \eqref{eq_def_functional} satisfies the conditions (i)-(iii), hence it is an
associated affine invariant of a
definition of
volume. 
\end{Lemma}

\proof
It is easy to check that $ \mathcal{V}$ is invariant under $\mathrm{GL}(V)$. Let us compute $ \mathcal{V}(B)$ for an
ellipsoid $B$. Since we already know that $\mathcal{V}$ is invariant under
$\mathrm{GL}(V)$, we may
choose a Euclidean scalar product and take $B$ as its unit ball. By Proposition \ref{prop_uniqueness}, the optimal body
$L$
in 
Theorem \ref{mainthm_convexity} is $\mathrm{SO}(n)$-invariant, hence a multiple of $B$, say $L=\lambda B$. The condition
on the
mixed volumes translates to $\lambda^{n-1} \omega_n=1$. The projection body operator is homogeneous of
degree $n-1$, and $\Pi B=\omega_{n-1}B^\circ$ hence 
\begin{displaymath}
  \mathcal{V}(B)=\frac{\omega_n^{n-1}}{\omega_{n-1}^n} \lambda^{n(n-1)} \svol(B \times \Pi
B)=\frac{\omega_n^{n-1}}{\omega_{n-1}^n} \frac{1}{\omega_n^n} \vol(B) \omega_{n-1}^n
\vol(B^\circ)=\omega_n.
\end{displaymath}

Next, suppose that $B \subset B'$, which implies that $B'^\circ \subset B^\circ$. Any probability
measure $\nu$ on $B'^\circ$ can be considered as a probability measure on $B^\circ$ and $\Gamma_\nu B^\circ=\Gamma_\nu
B'^\circ$.
Taking the maximum over such measures gives 
\begin{align*}
\frac{ \mathcal{V}(B)}{\vol B} & =\frac{\omega_n^{n-1}}{\omega_{n-1}^n} \left(\frac{n}{2}\right)^n \max_{\nu \in
\Prob(B^\circ)} \vol
\Gamma_\nu
B^\circ \\
& \geq \frac{\omega_n^{n-1}}{\omega_{n-1}^n}
\left(\frac{n}{2}\right)^n \max_{\nu \in \Prob(B'^\circ)} \vol \Gamma_\nu B'^\circ\\
& =\frac{ \mathcal{V}(B')}{\vol
B'}. 
\end{align*}
\endproof

In order to finish the proof of Theorem \ref{mainthm_convexity}, it remains to show that the corresponding definition
of volume is convex. 

\begin{Proposition}
 The definition of volume $\mu$ with associated affine functional $\mathcal{V}$ defined by \eqref{eq_def_functional} is
extendibly convex.
\end{Proposition}

\proof
Recall first that the definition of volume $\mu$ induces on each normed vector space $(V,B)$ a $k$-volume density
$\mu_k:\Lambda_s^k V \to \R$ by the formula 
\begin{displaymath}
 \mu_k(v_1 \wedge \ldots \wedge v_k):=\mu_W(v_1 \wedge \ldots \wedge v_k),
\end{displaymath}
where $W
\subset V$ is the $k$-plane spanned by $v_1,\ldots,v_k$, endowed with the induced norm.  

Using Theorem \ref{mainthm_gamma_mu}, we obtain the following explicit formula for $\mu_k$:
\begin{multline*}
 \mu_k(v_1 \wedge \ldots \wedge v_k)=\frac{\omega_k^{k-1}}{\omega_{k-1}^k} \frac{k^k}{k!} \max_{\nu \in \Prob((W \cap
B)^\circ)}\\ 
\left\{\int_{(W \cap B)^\circ}
\cdots \int_{(W \cap B)^\circ} \left|
\langle \eta_1 \wedge \ldots \wedge \eta_k,v_1 \wedge \ldots \wedge v_k\rangle\right| d\nu(\eta_1) \ldots
d\nu(\eta_k)\right\}.
\end{multline*}

We define a function $\tilde \mu_k$ on $\Lambda^kV$ by 
\begin{equation} \label{eq_tilde_mu}
\tilde \mu_k(\tau):=\frac{\omega_k^{k-1}}{\omega_{k-1}^k} \frac{k^k}{k!} \max_{\nu \in
\Prob(B^\circ)}\left\{\int_{B^\circ} \cdots \int_{B^\circ} \left|
\langle \xi_1
\wedge
\ldots \wedge \xi_k,\tau\rangle\right| d\nu(\xi_1) \ldots d\nu(\xi_k)\right\}, 
\end{equation}
where $\tau \in \Lambda^kV.$ Clearly it is convex. It remains to show that the restriction of $\tilde \mu_k$ to
the Grassmann cone of simple $k$-vectors equals $\mu_k$.  

Let $0 \neq \tau:=v_1 \wedge \ldots \wedge v_k \in \Lambda^k_sV$ and $W:=\langle v_1,\ldots,v_k\rangle$. Let $\iota:W
\to V$ be the inclusion. The dual map $\iota^*:V^* \to W^*$ is onto. Let $B':=W \cap B \subset W$ and $B'^\circ \subset W^*$ its polar. Then 
\begin{displaymath}
B'^\circ = \iota^*(B^\circ). 
\end{displaymath}
We may consider $\tau$ as an element of $\Lambda^kW$. Then
\begin{align*}
\tilde \mu_k & (\tau) = \frac{\omega_k^{k-1}}{\omega_{k-1}^k} \frac{k^k}{k!} \max_{\nu \in
\Prob(B^\circ)}\left\{\int_{B^\circ} \cdots \int_{B^\circ}
\left| \langle
\xi_1 \wedge \ldots \wedge \xi_k,\iota_*(\tau)\rangle\right| d\nu(\xi_1) \ldots d\nu(\xi_k)\right\} \\
& = \frac{\omega_k^{k-1}}{\omega_{k-1}^k} \frac{k^k}{k!} \max_{\nu \in \Prob(B^\circ)}\left\{\int_{B^\circ} \cdots
\int_{B^\circ} \left| \langle \iota^*
\xi_1 \wedge
\ldots \wedge \iota^* \xi_k,\tau\rangle\right| d\nu(\xi_1) \ldots d\nu(\xi_k)\right\} \\
& = \frac{\omega_k^{k-1}}{\omega_{k-1}^k} \frac{k^k}{k!} \max_{\nu \in \Prob(B^\circ)}\left\{\int_{B'^\circ}
\cdots \int_{B'^\circ} \left|
\langle \eta_1
\wedge \ldots \wedge \eta_k,\tau\rangle\right| d(\iota^*)_*\nu(\eta_1) \ldots d(\iota^*)_*\nu(\eta_k)\right\} \\
& = \frac{\omega_k^{k-1}}{\omega_{k-1}^k} \frac{k^k}{k!} \max_{\nu \in \Prob(B'^\circ)}\left\{\int_{B'^\circ}
\cdots \int_{B'^\circ } \left|
\langle \eta_1 \wedge \ldots \wedge \eta_k,\tau\rangle\right| d\nu(\eta_1) \ldots d\nu(\eta_k)\right\}\\
& = \mu_k(\tau),
\end{align*}
where the equality in the second to last line follows from the fact that the push-forward map $(\iota^*)_*:\Prob(B^\circ)
\to \Prob(B'^\circ)$ is onto.
\endproof

\section{The isoperimetrix}

In this section, we will describe the isoperimetrix for the new definition of
volume $\mu$. Let us first
recall the definition and construction of the isoperimetrix in general, referring to \cite{alvarez_thompson,
thompson_book96}
for more information. 

Let $\mu$ be a definition of volume, $(V,B)$ a normed vector space of dimension $n$ and suppose that the induced
$(n-1)$-volume density $\mu_{n-1}$ is convex. We
can integrate $\mu_{n-1}$ over $(n-1)$-dimensional submanifolds in $V$, in particular over the boundary of a compact
convex body $K$ (this makes sense even if $\partial K$ is not smooth). In this way we obtain the
surface
area $A_\mu(K)$ with respect to $\mu$. The isoperimetrix $\mathbb{I}_\mu B$ is the unique centrally
symmetric compact convex body in $V$ such that 
\begin{displaymath}
 A_\mu(K)=n V_1(K,\mathbb{I}_\mu B)
\end{displaymath}
for all $K$. 

Let us recall the construction of the isoperimetrix. The function $\mu_{n-1}:\Lambda^{n-1}V \to \R$ is convex and
$1$-homogeneous by assumption. The volume form on $V$ induces an isomorphism $\Lambda^{n-1}V \cong V^*$. We thus get a
convex and $1$-homogeneous function on $V^*$, which is the support function of the isoperimetrix. 

In the case of Busemann's definition of volume, the isoperimetrix is (up to a scalar) the polar of the
intersection body of $B$. The isoperimetrix of the Holmes-Thompson definition of volume is (again up to a scalar) the
projection body of the polar of the unit ball.

For a finite positive measure $\tau$ on a vector space $V$, we define the convex body $\Gamma \tau \subset V$ by 
\begin{displaymath}
 h(\Gamma \tau, \xi):=\int_V |\langle \xi,u\rangle| d\tau(u)
\end{displaymath}
and call it the {\it centroid body of $\tau$}. If $\tau$ is the normalized volume measure of a compact convex body $K$,
then
$\Gamma \tau=\Gamma K$ is the usual centroid body of $K$. 

\begin{Proposition}
 The isoperimetrix of the definition of volume $\mu$ from Theorem \ref{mainthm_convexity} is given by the following: for
each probability measure $\nu$ on $V^*$, let $\nu^\# \in \Prob(V)$ be the push-forward of $\nu \times \cdots \times \nu$
under
the map
\begin{align*}
 \underbrace{V^* \times \cdots \times V^*}_{n-1} & \to \Lambda^{n-1} V^* \cong V,\\
 \xi_1,\ldots,\xi_{n-1} & \mapsto \xi_1 \wedge \ldots \wedge \xi_{n-1}.
\end{align*}
Then the isoperimetrix of $B$ with respect to $\mu$ is given by 
\begin{displaymath}
 \mathbb{I}_\mu B = c_n \conv \bigcup_{\nu \in \Prob(B^\circ)} \Gamma \nu^\#,
\end{displaymath}
where 
\begin{displaymath}
 c_n:=\frac{\omega_{n-1}^{n-2}}{\omega_{n-2}^{n-1}} \frac{(n-1)^{n-1}}{(n-1)!}.
\end{displaymath}

\end{Proposition}

\proof
By \eqref{eq_tilde_mu}, $\mu_{n-1}=\tilde
\mu_{n-1}$ is $c_n$ times the maximum of the support functions of the $\Gamma \hat \nu$, where $\hat \nu$
is
the push-forward of $\nu \in \Prob(B^\circ)$ under the map $(V^*)^{n-1} \to \Lambda^{n-1}V^*$.  

Using the identification $\Lambda^{n-1}V^* \cong V$, we get that the support function of $\mathbb{I}_\mu B$ is
$c_n$ times the maximum of the support functions of the $\Gamma \nu^\#$, where $\nu \in \Prob(B^\circ)$. The
proof
is
finished by noting that the support function of the convex hull of a union of compact convex sets is the supremum of the
support functions. 
\endproof

\section{Affine inequalities}

\begin{Proposition} \label{prop_zonoid_case}
 Let $A$ a be a compact convex body. Then 
\begin{displaymath}
 \mathcal{V}(\Pi A) \leq \frac{\omega_n^{n-1}}{\omega_{n-1}^n} \vol \Pi A (\vol A)^{1-n}.
\end{displaymath}
Equality holds if and only if $A$ is homothetic to a projection body.  
\end{Proposition}

\proof
For each $L$  with $V_1(L,\Pi A)=1$ we find, using a well-known symmetry property of the projection
body operator (\cite{lutwak86}, Lemma 6)
\begin{displaymath}
\vol(\Pi L) \vol(A)^{n-1} \leq V_1(A,\Pi L)^n =V_1(L,\Pi A)^n= 1. 
\end{displaymath}
Equality holds if and only if $\Pi L$ and $A$ are homothetic. Taking the supremum (which is actually a maximum by
Theorem \ref{mainthm_gamma_mu}) over all such $L$ gives 
\begin{displaymath}
 \mathcal{V}(\Pi A) \leq \frac{\omega_n^{n-1}}{\omega_{n-1}^n} \vol(\Pi A)(\vol A)^{1-n}
\end{displaymath}
with equality if and only if $A$ is homothetic to a projection body.
\endproof

\begin{Corollary}
If $n=2$, then 
\begin{displaymath}
  \mathcal{V}(B)=\omega_2=\pi
\end{displaymath}
for all unit balls $B$. In particular, the $2$-volume density induced by $\mathcal{V}$ is Busemann's $2$-density. 
\end{Corollary}

\proof
In the two-dimensional case,  every centrally symmetric body is the projection body of some compact convex body.
We may thus write $B=\Pi A$ with $A$ centrally
symmetric. Since $\Pi A=2 JA$, where $J$ is rotation by $\frac{\pi}{2}$, it
follows that 
\begin{displaymath}
 \mathcal{V}(B)=\frac{\pi}{4} \frac{\vol(2 JA)}{\vol A}=\pi.
\end{displaymath}
\endproof

\begin{Proposition}
For all unit balls $B$,
\begin{displaymath}
  \mathcal{V}(B) \geq  \mathcal{V}^{ht}(B)
\end{displaymath}
with equality precisely for ellipsoids.
\end{Proposition}

\proof
Recall that the curvature image of $B$ is the unique (up to translations) compact convex body $\Lambda B$ with surface
area measure 
\begin{displaymath}
 dS_{n-1}(\Lambda B,\cdot)=h(B,\cdot)^{-n-1}d\sigma,
\end{displaymath}
where $\sigma$ is the spherical Lebesgue measure.

Let 
\begin{displaymath}
 L:=\frac{1}{(\vol B)^\frac{1}{n-1}} \Lambda B. 
\end{displaymath}
Using (\cite{lutwak86}, Lemmas 3 and 5),  
\begin{displaymath}
 V_1(L,B)=\frac{1}{\vol B} V_1(\Lambda B,B)=1
\end{displaymath}
and 
\begin{displaymath}
 \Pi L = \frac{1}{\vol B} \Pi \Lambda B=\frac{n+1}{2} \Gamma B^\circ.
\end{displaymath}

Using the Busemann-Petty centroid inequality \eqref{eq_busemann_petty_centroid} we get 
\begin{displaymath}
 \mathcal{V}(B) \geq \frac{\omega_n^{n-1}}{\omega_{n-1}^n} \svol(B \times \Pi L)  \geq  \frac{1}{\omega_n} \svol(B
\times B^\circ) =
\mathcal{V}^{ht}(B).
\end{displaymath}
\endproof

\begin{Proposition} \label{prop_equivalent_conj}
The assertion $\mathcal{V}(B) \geq \mathcal{V}^b(B)=\omega_n$ for all $B$ is equivalent to Petty's conjectured
projection inequality \ref{conj_petty}. 
\end{Proposition}

\proof
Set $L:=(\vol B)^{-\frac{1}{n-1}}B$. Then $V_1(L,B)=1$ and hence 
\begin{displaymath}
 \mathcal{V}(B) \geq \frac{\omega_n^{n-1}}{\omega_{n-1}^n} \svol(B \times \Pi L)=\frac{\omega_n^{n-1}}{\omega_{n-1}^n}
(\vol B)^{1-n} \vol \Pi B.
\end{displaymath}
Assuming Petty's conjectured projection inequality, the right hand side is bounded from below by $\omega_n$. 

Conversely, let $A$ be a compact convex body. Assuming $\mathcal{V}(\Pi A) \geq \omega_n$, 
Proposition \ref{prop_zonoid_case} implies that  
\begin{displaymath}
\frac{\omega_n^{n-1}}{\omega_{n-1}^n} \frac{\vol(\Pi A)}{(\vol A)^{n-1}} \geq  \mathcal{V}(\Pi A) \geq \omega_n,   
\end{displaymath}
with equality for ellipsoids. This is Petty's conjectured projection inequality.
\endproof


\end{document}